\documentclass[11pt,amssym,twoside]{article}
\usepackage{amssymb}
\newtheorem{thm}{Theorem}[section]
\newtheorem{mthm}[thm]{Main Theorem}

\newtheorem{conj}[thm]{Conjecture}
\newtheorem{prop}[thm]{Proposition}

\newtheorem{defn}[thm]{Definition}

\newtheorem{rem}[thm]{Remark}

\newcommand{\To}{\longrightarrow}
\newcommand{\Proj}{\mathbb P}
\newcommand{\Z}{\mathbb Z}

\newcommand{\F}{\mathbb F}
\newcommand{\D}{\mathbb D}
\newcommand{\Q}{\mathbb Q}
\begin{document}

\title{Deformation of Outer Representations of Galois Group}
\author{Arash Rastegar}
%\address{}%
%\email{rastegar@sharif.edu}%

%\thanks{}%
%\subjclass{}%
%\keywords{}%

%\date{}%
%\dedicatory{}
%\commby{}%
% ----------------------------------------------------------------
\maketitle
\begin{abstract}
To a hyperbolic smooth curve defined over a number-field one
naturally associates an "anabelian" representation of the
absolute Galois group of the base field landing in outer
automorphism group of the geometric fundamental group. In this
paper, we introduce several deformation problems for Lie-algebra
versions of the above representation and show that, this way we
get a richer structure than those coming from deformations of
"abelian" Galois representations induced by the Tate module of
associated Jacobian variety. To serve our purpose, we develop an
arithmetic deformation theory of graded Lie algebras with finite
dimensional graded components.
\end{abstract}
%------------------------------------------------------------------
\section*{Introduction}
One canonically associates to a smooth curve $X$ defined over a
number field $K$ a continuous group homomorphism
$$
\rho_X:Gal(\bar {K}/K) \To Out(\pi_1(X_{\bar {K}}))
$$
where $Out(\pi_1(X_{\bar {K}})) $ denotes the quotient of the
automorphism group $Aut(\pi_1(X_{\bar {K}})) $ by inner
automorphisms of the geometric fundamental group. By a conjecture
of Voevodsky and Matsumoto the outer Galois representation is
injective when topological fundamental group of $X$ is
nonabelian. Special cases of this conjecture are proved by Belyi
for $\mathbb {P}^1-\{0,1,\infty\}$ [Bel], by Voevodsky in cases of
genus zero and one [Voe], and by Matsumoto for affine $X$ using
Galois action on profinite braid groups [Mat1]. The importance of
the representation $\rho_X$ is due to Grothendieck's anabelian
conjectures [Gro1]. Mochizuki based on research by Nakamura and
Tamagawa proved Grothendieck's "Hom conjecture" which implies that
$\rho_X$ determines $X$ completely. More precisely, for $X$ and
$X'$ hyperbolic curves, the following natural map is a one-to-one
correspondence [Moc]:
$$
Isom_{K}(X,X')\To Out_{Gal(\bar {K}/K)}(Out(\pi_1(X_{\bar
{K}})),Out(\pi_1(X'_{\bar {K}})))
$$
Here $Out_{Gal(\bar {K}/K)}$ denotes the set of Galois
equivariant isomorphisms between the two profinite groups divided
by inner action of the second component.

The induced  pro-$l$ representation
$$
\rho^l_X:Gal(\bar {K}/K) \To Out(\pi^l_1(X_{\bar {K}}))
$$
after abelianization of the pro-$l$ fundamental group induces the
standard Galois representation associated to Tate module of the
Jacobian variety of $X$. Curves with abelian fundamental group are
not interesting here, because the outer representation does not
give any new information. After dividing $\pi^l_1(X_{\bar {K}})$
by its Frattini subgroup, or by mod-$l$ reduction of the
abelianized representation, one obtains a mod-$l$ representation
$$
\bar{\rho}^l_X:Gal(\bar {K}/K) \To GSp(2g,\F_l).
$$
We are interested in the space of deformations of the
representation $\rho^l_X$ fixing the mod-$l$ reduction
$\bar{\rho}^l_X$. In order to make sense of deforming a
representation landing in $Out(\pi^l_1(X_{\bar {K}}))$ we will
translate the outer representation of the Galois group to the
language of graded Lie-algebras.

The classical Schlessinger criteria for deformations of functors
on Artin local rings is used for deformation of the Galois action
on the abelianization of the pro-$l$ fundamental group which is
the same as etale cohomology. Using Schlessinger criteria, we will
construct universal deformation rings parameterizing all liftings
of the mod-$l$ representation $\bar{\rho}^l_X$ to actions of the
Galois group on graded Lie-algebras over $\Z_l$ with
finite-dimensional graded components. We also show that this
deformation theory is equivalent to deformation of abelian
representations of the Galois group.

To use the full power of outer representations, we deform the
corresponding Galois-Lie algebra representation to the graded Lie
algebra associated to weight filtration on outer automorphism
group of the pro-$l$ fundamental group. We construct a
deformation ring parameterizing all deformations fixing the
mod-$l$ Lie-algebra representation. Alongside, we develop an
arithmetic theory of deformations of Lie-algebras. The main point
we are trying to raise in this paper is that for a hyperbolic
curve $X$ the $l$-adic Lie-algebra representation we associate to
a hyperbolic curve $X$ contains more information than the
associated abelian $l$-adic representation.

The main obstacle in generalizing this method is the fact that
fundamental groups of curves are one relator groups and therefore
very similar to free groups. This makes it possible to mimic many
structures which work for free groups in the case of such
fundamental groups. This is heavily used in the course of
computations in this paper.

%------------------------------------------------------------------
\section{Background material}

The study of outer representations of the Galois group has two
origins. One root is the theme of anabelian geometry introduced by
Grothendieck [Gro1] which lead to results of Nakamura, Tamagawa
and Mochizuki who solved the problem in dimension one [Moc]. The
second theme which is originated by Deligne and Ihara
independently deals with Lie-algebras associated to the pro-$l$
outer representation [Del] [Iha]. This leads to a partial proof of
a conjecture by Deligne [Hai-Mat]. In the first part, we will
review the weight filtration introduced by Oda (after Deligne and
Ihara) and a circle of related results.

% ----------------------------------------------------------------
\subsection{Weight filtration on $\widetilde {Out}(\pi^l_1(X_{\bar {K}}))$}

By a fundamental result of Grothendieck [Gro2] the pro-$l$
geometric fundamental group $\pi^l_1(X_{\bar {K}})$ of a smooth
algebraic curve $X$ over $\bar {K}$ is isomorphic to the pro-$l$
completion of its topological fundamental group, after extending
the base field to the field of complex numbers. The topological
fundamental group of a Riemann surface of genus $g$ with $n$
punctured points has the following standard presentation:
$$
\Pi_{g,n}\cong <a_1,...,a_g,b_1,...,b_g,c_1,...,c_n>/
\prod_{i=1}^g [a_i,b_i]\prod_{j=1}^n c_j.
$$

Let $Aut(\pi^l_1(X_{\bar {K}}))$ denote the group of continuous
automorphisms of the pro-$l$ fundamental group of $X$ and
$Out(\pi^l_1(X_{\bar {K}}))$ denote its quotient by the subgroup
of inner automorphisms. We will induce filtrations on particular
subgroups of these two groups.

Let $\bar X$ denote the compactification of $X$ obtained by
adding finitely many points. $\bar X$ is still defined over $K$.
Let $\lambda:Aut(\pi^l_1(\bar {X}_{\bar {K}}))\to GL(2g,\Z_l)$
denote the map induced by abelianization. The natural actions of
$Aut(\pi^l_1(\bar X_{\bar {K}}))$ on cohomology groups
$H^i(\pi^l_1(\bar X_{\bar {K}}),\Z_l)$ are compatible with the
non-degenerate alternating form defined by the cup product:
$$
H^1(\pi^l_1(\bar X_{\bar {K}}),\Z_l)\times H^1(\pi^l_1(\bar
X_{\bar {K}}),\Z_l)\to H^2(\pi^l_1(\bar X_{\bar {K}}),\Z_l)\cong
\Z_l
$$
where the action of $\pi^l_1({\bar X}_{\bar {K}})$ on $\Z_l$ is
trivial. This shows that the image of $\lambda$ is contained in
$Sp(2g,\Z_l)$. One can prove that $\lambda$ is surjective and if
$\tilde {\lambda}$ denote the natural induced map
$$
\tilde{\lambda}:Out(\pi^l_1({\bar X}_{\bar {K}})) \longrightarrow
GSp(2g,\Z_l)
$$
there are explicit examples showing that the Galois representation
$\tilde{\lambda}\circ\rho^l_X$ does not fully determine the
original anabelian Galois representation [Asa-Kan].

Let $\widetilde {Aut}(\pi^l_1(X_{\bar {K}}))$ denote the braid
subgroup of $Aut(\pi^l_1(X_{\bar {K}}))$ which consists of those
elements taking each $c_i$ to a conjugate of a power
$c_i^{\sigma}$ for some $\sigma$ in $\Z_l^*$. There is a natural
surjective map
$$
\pi_X:\widetilde {Aut}(\pi^l_1(X_{\bar {K}})) \longrightarrow
Aut(\pi^l_1(\bar {X}_{\bar {K}}))
$$
Oda uses $\tilde {\lambda}$ to define and study natural
filtrations on $\widetilde{Aut}(\pi^l_1(X_{\bar {K}}))$ and its
quotient by inner automorphisms $\widetilde{Out}(\pi^l_1(X_{\bar
{K}}))$. In the special case of $X=\Proj^1-\{0,1,\infty\}$ this
is the same filtration as the filtration introduced by Deligne
and Ihara. This filtration is also used by Nakamura in bounding
Galois centralizers [Nak]. The central series of the pro-$l$
fundamental group induces a filtration on this group. This is not
the most appropriate filtration for non-compact $X$. In general,
we consider the weight filtration, namely the fastest decreasing
central filtration such that
$$
I^2 \pi^l_1(X_{\bar {K}})=<[\pi^l_1(X_{\bar {K}}),\pi^l_1(X_{\bar
{K}})],c_1,...,c_n>_{\textrm{norm}}
$$
where $<.>_{\textrm{norm}}$ means the closed normal subgroup
generated by these elements. For $m\geq 3$ we define
$$
I^m \pi^l_1(X_{\bar {K}})=<[I^i \pi^l_1(X_{\bar {K}}),I^j
\pi^l_1(X_{\bar {K}})]|i+j=m>_{\textrm{norm}}.
$$
The weight filtration induces a filtration on the automorphism
group of braid type by normal subgroups
$$
\widetilde {Aut}(\pi^l_1(X_{\bar {K}}))=I^0 \widetilde
{Aut}(\pi^l_1(X_{\bar {K}}))\supset I^1 \widetilde
{Aut}(\pi^l_1(X_{\bar {K}})) \supset ...\supset I^m \widetilde
{Aut}(\pi^l_1(X_{\bar {K}}))\supset ...
$$
$$
I^m \widetilde {Aut}(\pi^l_1(X_{\bar {K}}))=\{ \sigma \in
\widetilde {Aut}(\pi^l_1(X_{\bar {K}}))|x^{\sigma}x^{-1}\in
I^{m+1}\pi^l_1(X_{\bar {K}})\textrm{ for all }
x\in\pi^l_1(X_{\bar {K}})\}
$$
In fact $I^1 Aut(\pi^l_1(X_{\bar {K}}))$ coincides with the kernel
of $\lambda\circ \pi_X$.
\begin{prop} The weight filtration on $\widetilde
{Aut}(\pi^l_1(X_{\bar {K}}))$ satisfies
$$
[I^m \widetilde {Aut}(\pi^l_1(X_{\bar {K}})),I^n \widetilde
{Aut}(\pi^l_1(X_{\bar {K}}))]\subset I^{m+n}\widetilde
{Aut}(\pi^l_1(X_{\bar {K}}))
$$
for all $m$ and $n$, and induces a Lie-algebra structure on the
associated graded object $Gr^{\bullet}_I \widetilde
{Aut}(\pi^l_1(X_{\bar {K}}))=\bigoplus_{m\geq 1} gr^m \widetilde
{Aut}(\pi^l_1(X_{\bar {K}}))$. The graded pieces
$$
gr^m \widetilde {Aut}(\pi^l_1(X_{\bar {K}}))=I^m \widetilde
{Aut}(\pi^l_1(X_{\bar {K}}))/I^{m+1}\widetilde
{Aut}(\pi^l_1(X_{\bar {K}}))
$$
are free $\Z_l$-modules of finite rank for all positive $m$.
\end{prop}
The weight filtration on the automorphism group of pro-$l$
fundamental group induces a filtration on the outer automorphism
group of braid type
$$
\widetilde {Out}(\pi^l_1(X_{\bar {K}}))=I^0 \widetilde
{Out}(\pi^l_1(X_{\bar {K}}))\supset I^1 \widetilde
{Out}(\pi^l_1(X_{\bar {K}})) \supset ...\supset I^m \widetilde
{Out}(\pi^l_1(X_{\bar {K}}))\supset ...
$$
\begin{prop} The induced filtration on $\widetilde
{Out}(\pi^l_1(X_{\bar {K}}))$ satisfies
$$
[I^m \widetilde {Out}(\pi^l_1(X_{\bar {K}})),I^n \widetilde
{Out}(\pi^l_1(X_{\bar {K}}))]\subset I^{m+n}\widetilde
{Out}(\pi^l_1(X_{\bar {K}}))
$$
for all $m$ and $n$, and induces a Lie-algebra structure on the
associated graded object $Gr_I^{\bullet}\widetilde
{Out}(\pi^l_1(X_{\bar {K}}))=\bigoplus_{m\geq 1} gr^m \widetilde
{Out}(\pi^l_1(X_{\bar {K}}))$. The graded pieces
$$
gr^m \widetilde {Out}(\pi^l_1(X_{\bar {K}}))=I^m \widetilde
{Out}(\pi^l_1(X_{\bar {K}}))/I^{m+1}\widetilde
{Out}(\pi^l_1(X_{\bar {K}}))
$$
are finitely generated $\Z_l$-module for all positive $m$.
\end{prop}
These two propositions are proved in [kan].

%---------------------------------------------------------------
\subsection{Graded pieces of $\widetilde{Out}(\pi^l_1(X_{\bar {K}})$}

The explicit presentation of the fundamental group of a Riemann
surface given in the previous section implies that
$\pi^l_1(X_{\bar {K}})$ is a one-relator pro-$l$ group and
therefore very close to a free pro-$l$ group. The groups
$Aut(\pi^l_1(X_{\bar {K}}))$ and $Out(\pi^l_1(X_{\bar {K}}))$
also look very similar to automorphism group and outer
automorphism group of a free pro-$l$ group [Rib-Zal]. This can be
shown more precisely in the particular case of
$\widetilde{Out}(\pi^l_1(X_{\bar {K}}))$.

The graded pieces of $Gr^{\bullet}_I
\widetilde{Out}(\pi^l_1(X_{\bar {K}}))$ can be completely
determined in terms of the graded pieces of $Gr^{\bullet}_I
\pi^l_1(X_{\bar {K}})$ which are free $\Z_l$-modules. In fact,
$Gr^{\bullet}_I \pi^l_1(X_{\bar {K}})$ is a free Lie-algebra over
$\Z_l$ generated by images of $a_i$'s and $b_i$'s in $gr^1
\pi^l_1(X_{\bar {K}})$ for $1\leq i\leq g$ and $c_j$'s in $gr^2
\pi^l_1(X_{\bar {K}})$ for $1\leq j\leq n-1$. We denote these
generators by $\bar{a}_i$,$\bar{b}_i$ and $\bar{c}_j$
respectively.

Let $g_m$ denote the following injective $\Z_l$-linear
homomorphism
$$
g_m:gr^m \pi^l_1(X_{\bar {K}}) \longrightarrow (gr^{m+1}
\pi^l_1(X_{\bar {K}}))^{2g} \times (gr^m \pi^l_1(X_{\bar {K}}))^n
$$
$$
g\mapsto ([g,\bar{a}_i])_{1\leq i\leq g}\times
([g,\bar{b}_i])_{1\leq i\leq g}\times (g)_{1\leq j\leq n}
$$
and $f_m$ denote the following surjective $\Z_l$-linear
homomorphism
$$
f_m:(gr^{m+1} \pi^l_1(X_{\bar {K}}))^{2g} \times (gr^m
\pi^l_1(X_{\bar {K}}))^n \longrightarrow gr^{m+2} \pi^l_1(X_{\bar
{K}})
$$
$$
(r_i)_{1\leq i\leq g}\times (s_i)_{1\leq i\leq g}\times
(t_j)_{1\leq j\leq n}\mapsto \sum_{i=1}^g
([\bar{a}_i,s_i]+[r_i,\bar{b}_i])+\sum_{j=1}^n [t_j,\bar{c}_j].
$$
\begin{prop} The graded pieces of $Gr^{\bullet}_I
\widetilde{Out}(\pi^l_1(X_{\bar {K}}))$ fit into the following
short exact sequence of $\Z_l$-modules
$$
gr^m\widetilde{Out}(\pi^l_1(X_{\bar {K}})) \hookrightarrow
(gr^{m+1} \pi^l_1(X_{\bar {K}}))^{2g} \times (gr^m
\pi^l_1(X_{\bar {K}}))^n/gr^m \pi^l_1(X_{\bar {K}})
\twoheadrightarrow gr^{m+2} \pi^l_1(X_{\bar {K}})
$$
where embedding of $gr^m \pi^l_1(X_{\bar {K}})$ inside $(gr^{m+1}
\pi^l_1(X_{\bar {K}}))^{2g} \times (gr^m \pi^l_1(X_{\bar
{K}}))^n$ is defined by $g_m$ and the final surjection is induced
by $f_m$.
\end{prop}

This computational tool helps to work with the graded pieces of
$Gr^{\bullet}_I \widetilde{Out}(\pi^l_1(X_{\bar {K}}))$ as fluent
as the graded pieces of $Gr^{\bullet}_I \pi^l_1(X_{\bar {K}})$. In
particular, it enabled Kaneko to prove the following profinite
version of the Dehn-Nielson theorem [Kan]:
\begin{thm}(Kaneko)
Let $X$ be a smooth curve defined over a number-field $K$ and let
$Y$ denote an embedded curve in $X$ obtained by omitting finitely
many $K$-rational points. Then the natural map
$$
\widetilde{Out}(\pi^l_1(Y))\longrightarrow
\widetilde{Out}(\pi^l_1(X_{\bar {K}}))
$$
is a surjection.
\end{thm}
This can be easily proved by diagram chasing between the
corresponding short exact sequences for $Y$ and $X$. The above
exact sequence first appeared in the work of Ihara [Iha] and then
generalized by Asada and Kaneko [Asa-Kan] [Kan].

%---------------------------------------------------------------
\subsection{Filtrations on the Galois group}

If all of the points in the complement $\bar {X}-X$ are
$K$-rational, then the pro-$l$ outer representation of the Galois
group lands in the braid type outer automorphism group
$$
\tilde{\rho}^l_X:Gal(\bar {K}/K) \longrightarrow
\widetilde{Out}(\pi^l_1(X_{\bar {K}}))
$$
and the weight filtration on the pro-$l$ outer automorphism group
induce a filtration on the absolute Galois group mapping to
$\widetilde{Out}(\pi^l_1(X_{\bar {K}}))$ and also an injection
between associated Lie algebras over $\Z_l$ defined by each of
these filtrations
$$
Gr^{\bullet}_{X,l}Gal(\bar {K}/K)\hookrightarrow Gr_I^{\bullet}
\widetilde{Out}(\pi^l_1(X_{\bar {K}})).
$$

\begin{prop}
Let $X$ and $Y$ denote smooth curves over $K$ and let $\phi:X\to
Y$ denote a morphism also defined over $K$. Then $\phi$ induces a
commutative diagram of Lie algebras
$$
\begin{array}{ccc}
Gr^{\bullet}_{X,l}Gal(\bar {K}/K)& \hookrightarrow &
Gr_I^{\bullet} \widetilde{Out}(\pi^l_1(X_{\bar {K}})) \\
\downarrow \phi_{**}& & \downarrow \phi_* \\
Gr^{\bullet}_{Y,l}Gal(\bar {K}/K) & \hookrightarrow &
Gr_I^{\bullet} \widetilde{Out}(\pi^l_1(Y))
\end{array}
$$
If $\phi$ induces a surjection on topological fundamental groups,
then $\phi_*$ and $\phi_{**}$ will also be surjective.
\end{prop}
\textbf{Proof} The claim is true because $\phi_*$ is Galois
equivariant and graded pieces of $Gr_I^{\bullet}
\widetilde{Out}(\pi^l_1(X_{\bar {K}}))$ can be represented in
terms of exact sequences on graded pieces of $Gr_I^{\bullet}
\pi^l_1(X_{\bar {K}})$ [Kan]. $\square$
\begin{prop}
Let $X$ be an affine smooth curve $X$ over $K$ whose complement
has a $K$-rational point. Then there is a natural morphism
$$
Gr^{\bullet}_{X,l}Gal(\bar {K}/K)\to Gr^{\bullet}_{\Proj^1-\{
0,1,\infty\},l}Gal(\bar {K}/K).
$$
\end{prop}
\textbf{Proof} This is a consequence of theorem 3.1 in
[Mat1].$\square$
\begin{prop} There exists
a finite set of primes $S$ such that we have an isomorphism
$$
Gr^{\bullet}_{X,l}Gal(\bar {K}/K) \cong
Gr^{\bullet}_{X,l}Gal(K^{un}_S/K)
$$
where $Gal(K^{un}_S/K)$ denotes the Galois group of the maximal
algebraic extension unramified outside $S$
\end{prop}
\textbf{Proof.} Indeed, Grothendieck proved that the
representation $\tilde{\rho}^l_X$ factors through
$Gal(K^{un}_S/K)$ for a finite set of primes $S$. $S$ can be
taken to be primes of bad reduction of $X$ and primes over $l$
[Gro2]. This is also proved independently by Ihara in the special
case of $X=\Proj^1-\{ 0,1,\infty\}$ [Iha].$\square$

The importance of this result of Grothensieck is the fact that
$Gal(K^{un}_S/K)$ is a finitely generated profinite group
[Neu-Sch-Wei] and therefore, the moduli of its representations is
a scheme of finite type.

For a Lie algebra $L$ over $\Z_l$ let $Der(L)$ denote the set of
derivations, which are defined to be $\Z_l$-linear homomorphisms
$D:L\to L$ with
$$
D([u,v])=[D(u),v]+[u,D(v)]
$$
for all $u$ and $v$ in $L$, and let $Inn(L)$ denote the set of
inner derivations, which are defined to be derivations with
$D(u)=[u,v]$ for some fixed $v\in L$. Then we have the following
Lie algebra version of the outer representation of the Galois
group
$$
Gr^{\bullet}_{X,l}Gal(\bar {K}/K)\to Der(Gr^{\bullet}_I
\pi^l_1(X_{\bar {K}}))/Inn(Gr^{\bullet}_I \pi^l_1(X_{\bar {K}}))
$$
$$
\sigma\in gr^m Gal(\bar {K}/K) \mapsto (u\mapsto
\tilde{\sigma}(\tilde{u}).\tilde{u}^{-1})\textrm{ mod }
I^{m+n+1}\pi^l_1(X_{\bar {K}})
$$
where $u\in gr^n \pi^l_1(X_{\bar {K}})$  with $\tilde{u}\in I^n
\pi^l_1(X_{\bar {K}})$ and $\tilde{\sigma}\in I^m Gal(\bar
{K}/K)$ is a lift of $\sigma$. In fact, for a free graded
algebra, one can naturally associate a grading on the the algebra
of derivations. Let $L=\oplus_i L^i$ be a free graded Lie algebra
and let $D$ denote the derivation algebra of $L$. Then define
$$
D^i=\{ d\in D| D(L^j)\subset L^{i+j} \}.
$$
Then every element $d\in D$ is uniquely represented in the form
$d=\sum_i d^i$ with $d^i\in D^i$ such that for any $f\in L$ the
component $d^i f$ vanishes for almost all $i$. One can prove that
$$
[D^i,D^j]\subset D^{i+j}\textrm{ and }[D_1,D_2]^k=\sum_{i+j=k}
[D_1^i,D_2^j].
$$
One can mimic the same construction on the graded algebra
associated to $\pi^l_1(X_{\bar {K}})$ to get a graded algebra of
derivations [Tsu].

One shall notice that in case $X=\Proj^1-\{ 0,1,\infty\}$ the
group $\pi^l_1(X_{\bar {K}})$ is the pro-$l$ completion of a free
group with two generators. Ihara proves that the associated Lie
algebra $Gr^{\bullet}_I(\pi^l_1(X_{\bar {K}}))$ is also free over
two generators say $x$ and $y$ [Iha]. Now for $f$ in the $m$-th
piece of the grading of the Lie algebra, there is a unique
derivation $D_f\in Der (Gr^{\bullet}_I(\pi^l_1(X_{\bar {K}})))$
which satisfies $D_f(x)=0$ and $D_f(y)=[y,f]$. One can show that
$D_f(y)$ is non-zero for non-zero $m$ and that for any $\sigma\in
Gr_I^{\bullet}(\pi^l_1(X_{\bar {K}}))$ there exists a unique $f\in
Gr^{\bullet}_I(\pi^l_1(X_{\bar {K}}))$ with image of $\sigma$
being equal to $D_f$ [Mat2]. Now it is enough to let
$\sigma=\sigma_m$ the Soule elements, to get a non-zero image
$D_{f_{\sigma_m}}$.

\begin{conj} (Deligne) The graded Lie algebra
$(Gr^{\bullet}_{\Proj^1-\{ 0,1,\infty\},l}Gal(\bar
{\Q}/\Q))\otimes \Q_l$ is a free graded Lie algebra over $\Q_l$
which is generated by Soule elements and the Lie algebra
structure is induced from a Lie algebra over $\Z$ independent of
$l$. \end{conj}
\begin{rem}
It is reasonable to expect freeness to hold for $(Gr^{\bullet}_{X
,l}Gal(\bar {\Q}/\Q))\otimes \Q_l$.
\end{rem}
This implies that the above graded Lie algebra representation is
also injective. Ihara showed that Soule elements do generate
$Gr_{X,l}^{\bullet}(Gal(\bar {\Q}/\Q))\otimes \Q_l$ if one assumes
freeness of this Lie algebra [Iha]. Hain and Matsumoto proved the
same result without assuming any part of Deligne's conjecture
[Hai-Mat].

%-----------------------------------------------------------------
% ----------------------------------------------------------------
\section{Deformation theory}

We search for deformations of Galois outer representations which
are equipped with extra deformation data, for example with fixed
mod-$l$ representation. One can deform such representations both
in geometric and algebraic senses. A geometric example would be
given by $p$-adic deformation of a smooth curve $X$ defined over
a local field $K_{\wp}$. M. Kisin has considered such deformation
problems and has proved rigidity results in the mixed
characteristic case [Kis]. We are interested in algebraic
deformations which also work over global fields. In this part, we
will introduce several deformation problems for representations
landing in graded Lie algebras over $\Z_l$. We also develop an
arithmetic deformation theory for graded Lie algebras and
morphisms between them. In some cases universal deformations
exist and in some others, we are only able to construct a hull
which parameterizes all deformations.

%------------------------------------------------------------------
\subsection{Several deformation problems}

Let $k$ be a finite field of characteristic $p$ and let $\Lambda$
be any complete Noetherian local ring with residue field $k$. For
example $\Lambda$ can be $W(k)$, the ring of Witt vectors of $k$,
or $O$, the ring of integers of any local field with residue
field $k$. Let $C$ denote the category of Artinian local
$\Lambda$-algebras with residue field $k$. A covariant functor
from $C$ to $Sets$ is called pro-representable if it has the form
$$
F(A)\cong Hom_{\Lambda}(R_{univ},A)\hspace {0.5 in} A\in C
$$
where $R_{univ}$ is a noetherian complete local $\Lambda$-algebra
with maximal ideal $m_{univ}$ such that $R_{univ}/m_{univ}^n$ is
in $C$ for all $n$. There are a number of deformation functors
related to our problem.

The first deformation problem coming to mind is deforming actions
of Galois group on graded Lie algebras. Indeed, Galois group maps
to $\widetilde{Out}(\pi^l_1(X_{\bar {K}}))$ which acts on itself
and therefore on $Gr^{\bullet}_I\widetilde{Out}(\pi^l_1(X_{\bar
{K}}))$ by conjugation. By results of Kaneko, this action is
completely determined by the associated abelian Galois
representation [Kan]. This means that, the corresponding
deformation problem is pro-representable and we get exactly the
same deformation ring as in the abelian case.

In order to get a more delicate deformation theory, we could
deform the Lie algebra representation
$$
Gr^{\bullet}_{X,l} Gal(\bar {K}/K)\to
Gr^{\bullet}_I\widetilde{Out}(\pi^l_1(X_{\bar {K}}))
$$
of the Galois graded Lie algebra, induced by the outer
representation of the Galois group. We could define $D_{\bar
{\rho}}(A)$ for any Artinian local algebra $A$ to be the set of
isomorphism classes of deformations of the mod-$l$ reduction of
the above representation $\bar {\rho}$ to $Gr^{\bullet}_I
\widetilde{Out}(\pi^l_1(X_{\bar {K}}))\otimes A$. One is
interested in some variant of $D_{\bar {\rho}}(A)$ being
pro-representable.

The following Lie algebra version of an outer representation of
the Galois group
$$
Gr^{\bullet}_{X,l} Gal(\bar {K}/K)\to
Der(Gr^{\bullet}_I\pi^l_1(X_{\bar {K}}))/Inn(Gr^{\bullet}_I
\pi^p_1(X))
$$
which was defined in the first part is another candidate which
could be deformed.

%-------------------------------------------------------------------
\subsection{Galois actions on graded Lie algebras}

The method of proving modularity results by finding isomorphisms
between Hecke algebras and universal deformation rings as
originated by Wiles [Wil], can be reformulated in the language of
Lie algebras. One can define a Hecke-Lie algebra and a canonical
graded representation of the Galois-Lie algebra to Hecke-Lie
algebra which contains all the information of modular Galois
representations.

Let us first reformulate the theory of Galois representations in
the language Lie algebas. We start with elliptic curves. To each
elliptic curve $E$ defined over $\Q$ which has a rational point,
one associates a Galois outer representation
$$
Gal(\bar {\Q}/\Q) \longrightarrow Out(\pi^l_1(E-\{ 0\})).
$$
By analogy to Shimura-Taniyama-Weil conjecture, we expect this
representation be encoded in the representations
$$
Gal(\bar {\Q}/\Q) \longrightarrow Out(\pi^l_1(Y_0(N)))
$$
associated to modular curves $Y_0(N)$ which have a model over
$\Q$. By $Y_0(N)$ we mean the non-compactified modular curve of
level $N$ which is given as the quotient of the upper half-plane
by the congruence subgroup $\Gamma_0(N)$ of $SL_2(\Z)$ consisting
of matrices which are upper triangular modulo $N$.

For any smooth curve $X$ defined over $\Q$ the outer automorphism
group of braid type acts on $Gr^{\bullet}_I \widetilde
{Out}(\pi^l_1(X_{\bar {K}}))$ by conjugation and therefore for
each $m$ we get a Galois representation
$$
Gal(\bar {\Q}/\Q) \longrightarrow Aut(gr^m \widetilde
{Out}(\pi^l_1(X_{\bar {K}}))).
$$
In grade zero, we recover the usual abelian Galois representation
and in higher grades on can canonically construct this
representation by the grade-zero standard representation. Indeed,
for each $m\geq 1$ the isomorphism in proposition 1.3 is
$\widetilde {Out}(\pi^l_1(X_{\bar {K}}))$-equivariant. From this
we can determine the representation from the inner action of
$\widetilde {Out}(\pi^l_1(X_{\bar {K}}))$ on $Gr^{\bullet}_I
\pi^l_1(X_{\bar {K}})$. This action is fully determined by the
grade-zero action [Kan]. Therefore, the Galois representations
$$
Gal(\bar {\Q}/\Q) \longrightarrow Aut(gr^m
\widetilde{Out}(\pi^l_1(X_{\bar {K}})))
$$
are all determined by the abelian Galois representation
associated to $X$ over $\Q$. Together with Shimura-Tanyama-Weil
conjecture proved by Wiles and his collaborators [Wil] [Tay-Wil]
[Bre-Con-Dia-Tay] we get the following
\begin{thm}
Let $E$ be an elliptic curve over $\Q$ together with a rational
point $0\in E$. For each $m$ the Galois representation
$$
Gal(\bar {\Q}/\Q)\longrightarrow Aut(gr^m
\widetilde{Out}(\pi^l_1(E-\{0\})))
$$
appear as direct summand of the Galois representation
$$
Gal(\bar{\Q}/\Q) \longrightarrow Aut(gr^m
\widetilde{Out}(\pi^l_1(Y_0(N))))
$$
for some integer $N$.
\end{thm}

% -----------------------------------------------------------------
\subsection{Deformations of local graded Lie algebras}

We are interested in deforming the coefficient ring of graded Lie
algebras over $\Z_l$ of the form $L=Gr^{\bullet}_I \widetilde
{Out}(\pi^l_1(X_{\bar {K}}))$ and then deforming representations
of Galois graded Lie algebra
$$
Gr^{\bullet}_{X,l} Gal(\bar {K}/K)\to Gr^{\bullet}_I \widetilde
{Out}(\pi^l_1(X_{\bar {K}})).
$$
One can reduce the coefficient ring $\Z_l$ modulo $l$ and get a
graded Lie algebra $\bar L$ over $\F_l$ and a representation
$$
Gr^{\bullet}_{X,l} Gal(\bar {K}/K)\to \bar {L} .
$$
We look for liftings of this representation which is landing in
$\bar L$ among representations landing in graded Lie algebras
over Artin local rings $A$ of the form $L=\oplus_i L^i$ where
$L^0\cong GSp(2g,A)$ and $L^i$ is a finitely generated $A$-module
for positive $i$. In this section, we are only concerned with
deformations of Lie algebras and leave deformation of their
representations for the next section.

For a local ring $A\in C$ with maximal ideal $m$, the set of
deformations of $\bar L$ to $A$ is denoted by $D_c(\bar L,A)$ and
is defined to be the set of isomorphism classes of graded Lie
algebras $L/A$ of the above form which reduce to $\bar L$ modulo
$m$. In this notation $c$ stands for coefficients, since we are
only deforming coefficients not the Lie algebra structure. The
functor $D(\bar{L})$ as defined above is not a pro-representable
functor. As we will see, there exists a "hull" for this functor
(Schlessinger's terminology [Sch]) parameterizing all possible
deformations. One can also deform the Lie-algebra structure of
graded Lie algebras. The idea of deforming the Lie structure of
Lie algebras has been extensively used by geometers. For example,
Fialowski studied this problem in double characteristic zero case
[Fia]. In this paper, we are interested in double characteristic
$(0,l)$-version.

We shall first review cohomology of Lie algebras with the adjoint
representation as coefficients. Let $L$ be a graded $\Z_l$-algebra
and let $C^q(L,L)$ denote the space of all skew-symmetric
$q$-linear forms on a Lie algebra $L$ with values in $L$. Define
the differential
$$
\delta :C^q(L,L)\longrightarrow C^{q+1}(L,L)
$$
where the action of $\delta$ on a skew-symmetric $q$-linear form
$\gamma$ is a skew-symmetric $(q+1)$-linear form which takes
$(l_1,...,l_{q+1})\in L^{q+1}$ to
$$
\sum (-1)^{s+t-1}\gamma ([l_s,l_t],l_1,...,\hat l_s,...,\hat
l_t,...,l_{q+1})+\sum [l_u,\gamma(l_1,...,\hat l_u,...,l_{q+1})]
$$
where the first sum is over $s$ and $t$ with $1\leq s < t \leq
q+1$ and the second sum is over $u$ with $1\leq u\leq q+1$. Then
$\delta^2=0$ and we can define $H^q(L,L)$ to be the cohomology of
the complex $\{C^q(L,L),\delta\}$. If we put $C^m=C^{m+1}(L,L)$
and $H^m=H^{m+1}(L,L)$, then there exists a natural bracket
operation which makes $C=\oplus C^m$ a differential graded
algebra and $H=\oplus H^m$ a graded Lie algebra (see 1.2 of
[Fia-Fuc]). If $L$ is $\Z$-graded, $L=\oplus L(m)$, we say
$\phi\in C^q(L,L)(m)$ if for $l_i\in L(g_i)$ we have
$\phi(l_1,...,l_q)\in L(g_1+...+g_q-m)$. Then, there exists a
grading induced on the Lie algebra cohomology $H^q(L,L)=\oplus
H^q(L,L)(m)$.

Now we make an assumption for further constructions. Assume
$H^2(L,L)(m)$ is finite dimensional for all $m$ and consider the
algebra $\D_1=\Z_l\oplus\bigoplus_m H^2(L,L)(m)'$ where $'$ means
the dual over $\Z_l$. Fix a graded homomorphism of degree zero
$$
\mu :H^2(L,L)\longrightarrow C^2(L,L)
$$
which takes any cohomology class to a cocycle representing this
class. Now define a Lie algebra structure on
$$
\D_1\otimes L=L\oplus Hom^0(H^2(L,L),L)
$$
where $Hom^0$ means degree zero graded homomorphisms, by the
following bracket
$$
[(l_1,\phi_1),(l_2,\phi_2)]:=([l_1,l_2],\psi )
$$
where $\psi(\alpha )=\mu(\alpha)(l_1,l_2)+[\phi_1(\alpha),l_2]
+[\phi_2(\alpha),l_1]$. The Jacobi identity is implied by
$\delta\mu(\alpha)=0$. It is clear that this in an infinitesimal
deformation of $L$ and it can be shown that, up to an isomorphism,
this deformation does not depend on the choice of $\mu$. We shall
denote this deformation by $\eta_L$ after Fialowski and Fuchs
[Fia-Fuc].
\begin{prop} Any infinitesimal deformation of $\bar {L}$ to a
finite dimensional local ring $A$ is induced by pushing forward
$\eta_L$ by a unique morphism
$$
\phi:\Z_l\oplus\bigoplus_m H^2(L,L)(m)'\longrightarrow A.
$$
\end{prop}
\textbf {Proof.} This is the double characteristic version of
proposition 1.8 in [Fia-Fuc]. $\square$

Note that, in our case $H^2(L,L)$ is not finite dimensional. This
is why we restrict our deformations to the space of graded
deformations. Since $H^2(L,L)(0)$ is the tangent space of the
space of graded deformations, and the grade zero piece
$H^2(L,L)(0)$ is finite dimensional, the following version is more
appropriate:
\begin{prop} Any infinitesimal graded deformation of $\bar {L}$ to a
finite dimensional local ring $A$ is induced by pushing forward
$\eta^0_L$ by a unique morphism
$$
\phi:\Z_l\oplus H^2(L,L)(0)'\longrightarrow A.
$$
where $\eta^0_L$ denotes the restriction of $\eta_L$ to
$\Z_l\oplus H^2(L,L)(0)'$.
\end{prop}

Let $A$ be a small extension of $\mathbb {F}_l$ and $L$ be a
graded deformation of $\bar{L}$ over the base $A$. The
deformation space $D(\bar {L},A)$ can be identified with
$H^2(L,L)(0)$ which is finite dimensional. Therefore, by
Schlessinger criteria, in the subcategory $C'$ of $C$ consisting
of local algebras with $m^2=0$ for the maximal ideal $m$, the
functor $D(\bar {L},A)$ is pro-representable. This means that
there exists a unique map
$$
\Z_l\oplus H^2(L,L)(0)'\to R_{univ}
$$
inducing the universal infinitesimal graded deformation.

%-----------------------------------------------------------------
\subsection{Obstructions to deformations}

The computational tool used by geometers to study deformations of
Lie algebras is a cohomology theory of $K$-algebras where $K$ is
a field, which is developed by Harrison [Har]. This cohomology
theory is generalized by Barr to algebras over general rings
[Bar]. Here we use a more modern version of the latter introduced
independently by Andre and Quillen which works for general
algebras [Qui].

Let $A$ be a commutative algebra with identity over $\Z_l$ or any
ring $R$ and let $M$ be an $A$-module. By an $n$-long singular
extension of $A$ by $M$ we mean an exact sequence of $A$-modules
$$
0\to M\to M_{n-1}\to ...\to M_1 \to T \to R \to 0
$$
where $T$ is a commutative $\Z_l$-algebra and the final map a
morphism of $\Z_l$-algebras whose kernel has square zero. It is
trivial how to define morphisms and isomorphisms between $n$-long
singular extensions. Barr defines a group structure on these
isomorphism classes [Bar], which defines $H_{Barr}^n(A,M)$ for
$n>1$ and we put $H^1_{Barr}(A,M)=Der(A,M)$. Barr proves that for
a multiplicative subset $S$ of $R$ not containing zero
$$
H_{Barr}^n(A,M)\cong H_{Barr}^n(A_S,M)
$$
for all $n$ and any $A_S$-module $M$. According to this
isomorphism, the cohomology of the algebra $A$ over $\Z_l$ is the
same after tensoring $A$ with $\Q_l$ if $M$ is
$A\otimes\Q_l$-module. Thus one could assume that we are working
with an algebra over a field, and then direct definitions given by
Harrison would serve our computations better. Consider the complex
$$
0\rightarrow Hom(A,M) \rightarrow Hom(S^2A,M) \rightarrow
Hom(A\otimes A\otimes A,M)
$$
where $\psi\in Hom(A,M)$ goes to
$$
d_1\psi:(a,b)\mapsto a\psi(b)-\psi(ab)+b\psi(a)
$$
and $\phi\in Hom(S^2A,M)$ goes to
$$
d_2\phi(a,b,c)\mapsto a\phi(b,c)-\phi(ab,c)+\phi(a,bc)-c\phi(a,b)
$$
The cohomology of this complex defines $H^i_{Harr}(R,M)$ for
$i=1,2$. If $A$ is a local algebra with maximal ideal $m$ and
residue field $k$, the Harrison cohomology
$H^1_{Harr}(A,k)=(m/m^2)'$, which is the space of homomorphisms
$A\to k[t]/t^2$ such that $m$ is the kernel of the composition
$A\to k[t]/t^2\to k$.

Andre-Quillen cohomology is the same as Barr cohomology in low
dimensions and can be described directly in terms of derivations
and extensions. For any morphism of commutative rings $A \to B$
and $B$-module $M$ we denote the $B$-module of $A$-algebra
derivations of $B$ with values in $M$ by $Der_A (B,M)$. Let
$Ext^{inf}_A (B,M)$ denote the $B$-module of infinitesimal
$A$-algebra extensions of $B$ by $M$. The functors $Der$ and
$Ext^{inf}$ have transitivity property. Namely, given morphisms
of commutative rings $A\to B\to C$ and a $C$-module $M$, there is
an exact sequence
$$
0 \longrightarrow Der_B (C,M)\longrightarrow Der_A
(C,M)\longrightarrow Der_A (B,M) \hspace{1 in}$$
$$
\hspace{1 in}\longrightarrow Ext^{inf}_B (C,M)\longrightarrow
Ext^{inf}_A (C,M)\longrightarrow Ext^{inf}_A (B,M).
$$
The two functors $Der$ and $Ext^{inf}$ also satisfy flat
base-change property. Namely, given morphisms $A\to B$ and $A\to
A'$ if $Tor^A_1(A',B)=0$, then there are isomorphisms $Der_{A'}
(A'\otimes_A B,M)\cong Der_A (B,M)$ and $Ext^{inf}_{A'}
(A'\otimes_A B,M)\cong Ext^{inf}_A (B,M)$. Andre-Quillen
cohomology associates $Der_A (B,M)$ and $Ext^{inf}_A (B,M)$ to any
morphism of commutative rings $A \to B$ and $B$-module $M$ as the
first two cohomologies and extends it to higher dimensional
cohomologies such that transitivity and flat base-change extend
in the obvious way.

Let $A$ be an object in the category $C$ and $L\in Def(\bar
{L},A)$. The pair $(A,L)$ defines a morphism of functors $\theta
:Mor(A,B)\to Def(\bar {L},B)$. We say that $(A,L)$ is universal if
$\theta$ is an isomorphism for any choice of $B$. We say that
$(A,L)$ is miniversal if $\theta$ is always surjective, and gives
an isomorphism for $B=k[\varepsilon]/\varepsilon^2$. We intend to
construct a miniversal deformation of $\bar {L}$.

Consider a graded deformation with base in a local algebra $A$
with residue field $k=\F_l$. One can define a map
$$
\Phi_A:Ext^{inf}_{\Z_l}(A,k)\longrightarrow H^3(\bar {L},\bar
{L}).
$$
Indeed, choose an extension $0\to k\to B\to A \to 0$ corresponding
to an element in $Ext^{inf}_{\Z_l}(A,k)$. Consider the $B$-linear
skew-symmetric operation $\{ .,.\}$ on $\bar {L}\otimes_k B$
commuting with $[.,.]$ on $\bar {L}\otimes A$ defined by $\{
l,l_1\}=[l,\bar {l}_1]$ for $l$ in the kernel of $\bar
{L}\otimes_k B\to \bar {L}\otimes_k A$ which can be identified by
$\bar {L}$. Here $\bar {l}_l$ is the image of $l_1$ under the
projection map $B\to k$ tensored with $\bar L$ whose kernel is the
inverse image of the maximal ideal of $A$. The Jacobi expression
induces a multilinear skew-symmetric form on $\bar L$ which could
be regarded as a closed element in $C^3(\bar {L},\bar {L})$. The
image in $H^3(\bar {L},\bar {L})$ is independent of the choices
made.
\begin{thm}(Fialowski)
One can deform the Lie algebra structure on $L\otimes_k A$ to
$L\otimes_k B$ if and only if the image of the above extension
vanishes under the morphism $Ext^{inf}_{\Z_l}(A,k)\to H^3(\bar
{L},\bar {L})$.
\end{thm}
\textbf {Proof.} The proof presented in [Fia] and [Fia-Fuc] works
for algebras over fields of finite characteristic. $\square$

Using the above criteria for extending deformations, one can
follow the methods of Fialowski and Fuchs to introduce a
miniversal deformation for $\bar L$.
\begin{prop}
Given a local commutative algebra $A$ over $\Z_l$ there exists a
universal extension
$$
0 \longrightarrow Ext^{inf}_{\Z_l}(A,k)'\longrightarrow C
\longrightarrow A \longrightarrow 0
$$
among all extensions of $A$ with modules $M$ over $A$ with $mM=0$
where $m$ is the maximal ideal of $A$.
\end{prop}
This is proposition 2.6 in [Fia-Fuc]. Consider the canonical split
extension
$$
0\longrightarrow H^2(L,L)(0)'\longrightarrow \D_1\longrightarrow
k\longrightarrow 0.
$$
We will initiate an inductive construction of $\D_k$ such that
$$
0\longrightarrow Ext^{inf}_{\Z_l}(\D_k,k)'\longrightarrow \bar
{\D}_{k+1}\longrightarrow \D_k\longrightarrow 0.
$$
together with a deformation $\eta_k$ of $L$ to the base $\D_k$.
For $\eta_1$ take $\eta_L$, and assume $\D_i$ and $\eta_k$ is
constructed for $i\leq k$. Given a local commutative algebra
$\D_k$ with maximal ideal $m$, there exists a unique universal
extension for all extensions of $\D_k$ by $\D_k$-modules $M$ with
$mM=0$ of the following form
$$
0\longrightarrow Ext^{inf}_{\Z_l}(\D_k,k)'\longrightarrow C
\longrightarrow \D_k\longrightarrow 0.
$$
associated to the cocycle $f_k:S^2(\D_k)\to
Ext^{inf}_{\Z_l}(\D_k,k)'$ which is dual to the homomorphism
$$
\mu: Ext^{inf}_{\Z_l}(\D_k,k)\longrightarrow S^2(\D_k)'
$$
which takes a cohomology class to a cocycle from the same class.
The obstruction to extend $\eta_k$ lives in
$Ext^{inf}_{\Z_l}(\D_k,k)'\otimes H^3(\bar {L},\bar {L})$.
Consider the composition of the associated dual map
$$
\Phi'_k:H^3(\bar {L},\bar {L})'\longrightarrow
Ext^{inf}_{\Z_l}(\D_k,k)'.
$$
with $Ext^{inf}_{\Z_l}(\D_k,k)'\to  C$ and define $\D_{k+1}$ to
be the cokernel of this map. We get the following exact sequence
$$
0\longrightarrow (ker\Phi_k)'\longrightarrow \D_{k+1}
\longrightarrow \D_k\longrightarrow 0.
$$
We can extend $\eta_k$ to $\eta_{k+1}$. Now, taking a projective
limit of $\D_k$ we get a base and a formal deformation of $\bar
{L}$.
\begin{thm}
Let $\D$ denote the projective limit $\underline {\lim} \D_k$
which is a $\Z_l$-module. One can deform $\bar L$ uniquely to a
graded Lie algebra with base $\D$ which is miniveral among all
deformations of $\bar {L}$ to local algebras over $\Z_l$.
\end{thm}
\textbf {Proof.} This is the double characteristic version of
theorem 4.5 in [Fia-Fuc]. The same proof works here because
theorems 11 and 18 in [Har] which are used in the arguments of
Fialowski and Fuchs work for algebras over any perfect field.
$\square$

\begin{prop} (Fialowski-Fuchs) The base of the minversal
deformation of $\bar {L}$ is the zero locus of a formal map
$H^2(\bar {L},\bar {L})(0)\to H^3(\bar {L},\bar {L})(0)$.
\end{prop}
This is the graded version of proposition 7.2 in [Fia-Fuc].

% -----------------------------------------------------------------
\subsection{Deformations of graded Lie algebra representations}

In the previous section we discussed deformation theory of the
mod $l$ reduction of the Lie algebra $Gr^{\bullet}_I \widetilde
{Out}(\pi^l_1(X_{\bar {K}}))$. We shall mention the following
\begin{thm}  The cohomology
groups $H^i(Gr^{\bullet}_I \widetilde {Out}(\pi^l_1(X_{\bar
{K}})),Gr^{\bullet}_I \widetilde {Out}(\pi^l_1(X_{\bar
{K}})))(0)$ are finite dimensional for all non-negative integer
$i$.
\end{thm}
\textbf{Proof.} By a theorem of Labute $Gr^{\bullet}_I
\pi^l_1(X_{\bar {K}})$ is quotient of a finitely generated free
Lie algebra with finitely generated module of relations [Lab].
Therefore, the cohomology groups $H^i(Gr^{\bullet}_I
\pi^l_1(X_{\bar {K}}),Gr^{\bullet}_I \pi^l_1(X_{\bar {K}}))(0)$
are finite dimensional. Finite dimensionality of the cohomology of
$Gr^{\bullet}_I \widetilde {Out}(\pi^l_1(X_{\bar {K}}))$ follows
from proposition 1.3. $\square$

We are interested in deforming the following graded
representation of the Galois graded Lie algebra
$$
\rho:Gr^{\bullet}_{X,l} Gal(\bar {K}/K)\to Gr^{\bullet}_I
\widetilde {Out}(\pi^l_1(X_{\bar {K}}))
$$
among all graded representations which modulo $l$ reduce to the
graded representation
$$
\bar {\rho}:Gr^{\bullet}_{X,l} Gal(\bar {K}/K)\to \bar {L}
$$
where the Lie algebra $\bar L$ over $\F_l$ is the mod-$l$
reduction of $Gr^{\bullet}_I \widetilde {Out}(\pi^l_1(X_{\bar
{K}}))$. There are suggestions from the classical deformation
theory of Galois representations on how to get a representable
deformation functor [Til]. Let $D(\bar {\rho},A)$ denote the set
of isomorphism classes of Galois graded Lie algebra
representations to graded Lie algebras $L/A$ of the above form
which reduce to $\bar {\rho}$ modulo $m$. The first ingredient we
need to prove $D(\bar {\rho})$ is representable is that the
tangent space of the functor is finite-dimensional. The tangent
space of the deformation functor $D(\bar {\rho})$ for an object
$A\in C$ is canonically isomorphic to
$$
H^1(Gr^{\bullet}_{X,l} Gal(\bar {K}/K),Ad\circ \overline {\rho})
$$
where the Lie algebra module is given by the composition of $\bar
{\rho}$ with the adjoint representation of $Gr^{\bullet}_I
\widetilde {Out}(\pi^l_1(X/K))$. To get finite dimensionality, we
restrict ourselves to the graded deformations of the graded
representation $\bar {\rho}$.
\begin{thm}
$H^1(Gr^{\bullet}_{X,l} Gal(\bar {K}/K),Ad\circ \overline
{\rho})(0)$ is finite dimensional.
\end{thm}
\textbf {Proof.} The Galois-Lie representation $Gr^{\bullet}_I
Gal(\bar {K}/K)\to Gr^{\bullet}_I \widetilde
{Out}(\pi^l_1(X_{\bar {K}}))$ is an injection. Derivation
inducing cohomology commutes with inclusion of Lie algebras.
Therefore $H^1(Gr^{\bullet}_{X,l} Gal(\bar {K}/K),Ad\circ
\overline {\rho})(0)$ injects in $H^1(Gr^{\bullet}_I \widetilde
{Out}(\pi^l_1(X_{\bar {K}})),Ad\circ \overline {\rho})(0)$ which
is finite dimensional by previous theorem. $\square$

For a surjective mapping $A_1\to A_0$ of Artinian local rings in
$C$ such that the kernel $I\subset A_1$ satisfies $I.m_1=0$ and
given any deformation $\rho_0$ of $\bar {\rho}$ to $L\otimes A_0$
one can associate a canonical obstruction class in
$$
H^2(Gr^{\bullet}_{X,l} Gal(\bar {K}/K),I\otimes Ad\circ \bar
{\rho})\cong H^2(Gr^{\bullet}_{X,l} Gal(\bar {K}/K),Ad\circ \bar
{\rho})\otimes I
$$
which vanishes if and only if $\rho_0$ can be extended to a
deformation with coefficients $A_1$. Therefore, vanishing results
on second cohomology are important.
\begin{thm}
Suppose $Gr^{\bullet}_{X,l} Gal(\bar {K}/K)$ is a free Lie algera
over $\Z_l$, then the Galois cohomology $H^2(Gr^{\bullet}_{X,l}
Gal(\bar {K}/K),Ad\circ \bar {\rho})$ vanishes.
\end{thm}
\textbf {Proof.} The free Lie algebra $G=Gr^{\bullet}_{X,l}
Gal(\bar {K}/K)$ is rigid, and therefore has trivial infinitesimal
deformations. Thus, we get vanishing of its second cohomology:
$$
H^2(Gr^{\bullet}_{X,l} Gal(\bar {K}/K),Gr^{\bullet}_{X,l} Gal(\bar
{K}/K))=0.
$$
The injection of $G$ inside $L=Gr^{\bullet}_I \widetilde
{Out}(\pi^l_1(X_{\bar {K}}))$ as Lie-algebras over $\Z_l$ implies
that, the cohomology group $H^2(Gr^{\bullet}_{X,l} Gal(\bar
{K}/K),Ad\circ \rho)$ vanishes again by freeness of $G$. Let
$\bar {G}$ denote the reduction modulo $l$ of $G$ which is a free
Lie algebra over $\F_l$. The cohomology $H^2(G,\bar {G})$ is the
mod-$l$ reduction of $H^2(G,G)$, hence it also vanishes. So does
the cohomology $H^2(Gr^{\bullet}_{X,l} Gal(\bar {K}/K),Ad\circ
\bar {\rho})$ by similar reasoning.$\square$

We have obtained conceptual conditions implying $H_4$ of [Sch].
What we have proved can be summarized as follows.

\begin{mthm} Suppose that $Gr^{\bullet}_{X,l} Gal(\bar {K}/K)$ is a
free Lie algera over $\Z_l$. There exists a universal deformation
ring $R_{univ}=R(X,K,l)$ and a universal deformation of the
representation $\bar {\rho}$
$$
\rho^{univ}:Gr^{\bullet}_{X,l} Gal(\bar {K}/K)\longrightarrow
Gr^{\bullet}_I \widetilde {Out}(\pi^l_1(X_{\bar {K}}))\otimes
R_{univ}
$$
which is unique in the usual sense. If $Gr^{\bullet}_{X,l}
Gal(\bar {K}/K)$ is not free, then a mini-versal deformation
exists which is universal among infinitesimal deformations of
$\bar {\rho}$.
\end{mthm}
\begin{rem}
Note that, freeness of $Gr^{\bullet}_{X,l} Gal(\bar {K}/K)$ in the
special case of $K=\Q$ where filtration comes from punctured
projective curve $X=\mathbb {P}^1-\{ 0,1,\infty\}$ or a punctured
elliptic curve $X=E-\{ 0\}$ is implied by Deligne's conjecture.
\end{rem}

As in the classical case, the Lie algebra structure on
$Ad\circ\bar{\rho}$ induces a graded Lie algebra structure on the
cohomology $H^*(Gr^{\bullet}_{X,l} Gal(\bar {K}/K),Ad\circ \bar
{\rho})$ via cup-product, and in particular, a symmetric bilinear
pairing
$$
H^1(Gr^{\bullet}_{X,l} Gal(\bar {K}/K),Ad\circ \bar {\rho})\times
H^1(Gr^{\bullet}_{X,l} Gal(\bar {K}/K),Ad\circ \bar {\rho})
\hspace {1 in}
$$
$$
\hspace {3 in}\longrightarrow H^2(Gr^{\bullet}_{X,l} Gal(\bar
{K}/K),Ad\circ \bar {\rho})
$$
which gives the quadratic relations satisfied by the minimal set
of formal parameters of $R_{univ}/lR_{univ}$ for characteristic
$l$ different from $2$.
% -----------------------------------------------------------------
\subsection{Hecke-Teichmuller graded Lie algebras}

There is no general analogue for Hecke operators in the context
of Lie algebras constructed in this manner. What we need is an
analogue of Hecke algebra which contains all the information of
Galois outer representations associated to elliptic curves. In
fact, we will provide an algebra containing such information for
hyperbolic smooth curves of given topological type. From now on,
we assume that $2g-2+n>0$. The stack $M_{g,n}$ is defined as the
moduli stack of $n$-pointed genus $g$ curves.
\begin{defn}
A family of $n$-pointed genus $g$ curves over a scheme $S$ is a
proper smooth morphism $C\to S$ whose fibers are proper smooth
curves of genus $g$ together with $n$ sections $s_i:S\to C$ for
$i=1,...,n$ whose images do not intersect.
\end{defn}

The moduli stack is an algebraic stack over $Spec(\Z)$. One can
define the etale fundamental group of the stack $M_{g,n}$ in the
same manner one defines etale fundamental group of schemes. Oda
showed that the etale homotopy type of the algebraic stack
$M_{g,n}\otimes \bar {\Q}$ is the same as the analytic stack
$M^{an}_{g,n}$ and its algebraic fundamental group is isomorphic
to the completion $\widehat {\Gamma_{g,n}}$ of the Teichmuller
modular group, or the mapping class group of $n$-punctured genus
$g$ Riemann surfaces [Oda]:
$$
\pi_1(M_{g,n}\otimes \bar {\Q}) \cong \widehat {\Gamma_{g,n}}.
$$
Triviality of $\pi_2$ implies exactness of the following short
sequence for the universal family $C_{g,n}\To M_{g,n}$ over the
moduli stack
$$
0\To \pi_1(C,b)\To \pi_1(C_{g,n},b)\To \pi_1(M_{g,n},a) \To 0
$$
where $C$ is the fiber on $a$ and $b$ is a point on $C_{g,n}$.
Using this exact sequence, one defines the arithmetic universal
monodromy representation
$$
\rho_{g,n}:\pi_1(M_{g,n},a) \To Out(\pi_1(C))
$$
In fact, this is the completion of the natural map
$\Gamma_{g,n}\To Out(\Pi_{g,n})$. By composition with the natural
projection to outer automorphism group of the $l$-adic completion
$\Pi^l_{g,n}$ we get a representation
$$
\rho_{g,n}:\pi_1(M_{g,n},a) \To Out(\pi_1^l(C)).
$$
This map induces filtrations on $\pi_1(M_{g,n})$ and its subgroup
$\pi_1(M_{g,n}\otimes \bar{\Q})$ and an injection of $\Z_l$-Lie
algebras
$$
Gr^{\bullet}_I \pi_1(M_{g,n}\otimes \bar{\Q})\hookrightarrow
Gr^{\bullet}_I \pi_1(M_{g,n}).
$$
It is conjectured by Oda and proved by a series of papers by
Ihara, Matsumoto, Nakamura and Takao that the cokernel of the
above map after tensoring with $\Q_l$ is independent of $g$ and
$n$ [Iha-Nak] [Mat1] [Nak]. Note that $M_{0,3}\cong Spec(\Q)$.

On can think of $\pi_1(M_{g,n}\otimes \bar{\Q})$ as a replacement
for the Hecke algebra. This object has the information of all
outer representations associated to smooth curves over $\Q$.
Indeed, by fixing such a curve $C$ of genus $g\geq 2$ we have
introduced a point on the moduli stack $a\in M_{g,0}$ and thus a
Galois representation
$$
Gal(\bar{\Q}/\Q)\To \pi_1(M_{g,n})
$$
which splits the following short exact sequence
$$
0\To \pi_1(M_{g,n}\otimes \bar {\Q},a)\To \pi_1(M_{g,n}) \To
Gal(\bar {\Q}/\Q)\To 0.
$$
Combining with the arithmetic universal monodromy representation
we get
$$
Gal(\bar{\Q}/\Q)\To Out(\Pi^l_{g,n})
$$
which recovers the canonical outer representation associated to
$C$.
\begin{defn}
We define the Hecke-Teichmuller Lie algebra to be the image of the
following morphism of graded Lie algebras
$$
Gr^{\bullet}_I \pi_1(M_{g,n}) \To Gr^{\bullet}_I \widetilde
{Out}(\Pi^l_{g,n}).
$$
\end{defn}
We expect Hecke-Teichmuller Lie algebra to serve the role of
Hecke algebra in proving modularity results for elliptic curves
or other motivic objects.
%-------------------------------------------------------------------
\section*{Acknowledgements}
I would like to thank M. Kontsevich, O. Gabber, N. Nikolov, V.B.
Mehta, A.J. Parameswaran, D. Prasad, Y. Soibelman for enjoyable
conversations. Also I wish to thank Tata institute of fundamental
research and institute des hautes \`{e}tudes scientifiques for
their warm hospitality during which this work was prepared. I am
grateful to referees for comments and corrections.
%-------------------------------------------------------------------

Sharif University of Technology, e-mail: rastegar@sharif.ir


\begin{thebibliography}{999999}
%\bibliographystyle{amsplain}
%\bibliography{bibliography}

\bibitem[Asa-Kan]{Asa-Kan} M. Asada, M. Kaneko; {\em On the automorphism
group of some pro-$l$ fundamental groups,} Advanced studies in
Pure Mathematics Vol. 12,137-159(1987).

\bibitem[Bar]{Bar} M. Barr; {\em A cohomology theory for
commutative algebras I II,} Proc. A.M.S. 16,1379-1391(1965).

\bibitem[Bel]{Bel} G.V. Belyi; {\em On Galois extensions of a maximal
cyclotomic field,} Math USSR Izv. 14,247-256(1980).

\bibitem[Bre-Con-Dia-Tay]{Bre-Con-Dia-Tay} C. Breuil, B. Conrad, F. Diamond,
R. Taylor; {\em On the modularity of elliptic curves over $\Q$:
wilde 3-adic excercises,} J. amer. Math. Society 14 no.
4,843-939(2001).

\bibitem[Del]{Del} P. Deligne; {\em Le groupe fundamental de la droite
projective moins trois points,} in {\rm Galois groups over $\Q$}
79-297, edited by Y. Ihara, J.P. Serre, Springer 1989.

\bibitem[Fia]{Fia} A. Fialowski; {\em Deformations of Lie algebras,}
Math USSR Sbornik Vol. 55,467-473(1986).

\bibitem[Fia-Fuc]{Fia-Fuc} A. Fialowski, D. Fuchs; {\em Construction
of miniversal deformations of Lie algebras,} Funct. Anal.
161(1),76-110(1999).

\bibitem[Gro1]{Gro1} A. Grothendieck; {\em Esquisse d'un
programme,} manuscript, 1984.

\bibitem[Gro2]{Gro2} A. Grothendieck; {\em Revetement Etale et
Groupe Fondamental (SGA I),} LNM 224, Springer-Verlag 1971.

\bibitem[Hai-Mat]{Hai-Mat} R. Hain, M. Matsumoto; {\em Weighted
Completion of Galois Groups and Galois Actions on the Fundamental
Group of $\mathbb P^1 - \{ 0,1,\infty\}$} Compositio Math.
139,119-167(2003).

\bibitem[Har]{Har} D.K. Harrison; {\em Commutative algebras and
cohomology,} Trans. A.M.S. 104,191-204(1962).

\bibitem[Iha]{Iha} Y. Ihara; {\em Profinite braid groups, Galois
representations, and complex multiplications,} Ann. of Math.
123,43-106(1986).

\bibitem[Iha-Nak]{Iha-Nak} Y. Ihara, H. Nakamura; {\em On deformation
of mutually degenerate stable marked curves and Oda's problem,}
J. Reine Angenw. Math. 487,125-151 (1997).

\bibitem[Kan]{Kan} M. Kaneko; {\em Certain automorphism groups of
pro-$l$ fundamental groupsof punctured Riemann surfaces,} J. Fac.
Sci. Univ. Tokyo Sect. IA, Math. 36,363-372 (1989).

\bibitem[Kis]{Kis} M. Kisin; {\em Prime to $p$ fundamental groups
and tame Galois actions,} Ann. Inst. Fourier, Grenoble 50,
4,1099-1126(2000).

\bibitem[Lab]{Lab} J. Labute; {\em On the descending central
series of groups with single defining relation,} J. of Algebra,
14,16-23(1970).

\bibitem[Mat1]{Mat1} M. Matsumoto; {\em Galois representations on
profinite braid groups on curves,} J. Reine Angenw. Math.
474,169-219 (1996).

\bibitem[Mat2]{Mat2} M. Matsumoto; {\em On the Galois image in the
derivation algebra of $\pi_1$ of the projective line minus three
points,} Contemporary Math. 186,201-213(1995).

\bibitem[Moc]{Moc} S. Mochizuki; {\em The profinite Grothendieck
conjecture for closed hyperbolic curves over number fields,} J.
Math. Sci. Univ. Tokyo 3,571-627(1996).

\bibitem[Nak]{Nak} H. Nakamura; {\em Coupling of universal monodromy
representations of Galois Techmuller modular groups,} Math. Ann
304,99-119 (1996).

\bibitem[Neu-Sch-Win]{Neu-Sch-Win} J. Neukirch, A. Schmidt, Wingberg;
{\em Cohomology of number fields,} Springer 2000.

\bibitem[Oda]{Oda} T. Oda; {\em Etale homotopy type of the moduli
spaces of algebraic curves,} in "Geometric Galois Actions I"
London Math. Soc. Lect. Note Ser. 242,85-95 (1997).

\bibitem[Qui]{Qui} D. Quillen; {\em On the (co)-homology of
commutative rings,} in "Applications of commutative algebra"
Proc. Sympos. Pure. Math. Vol. XVII, Amer. Math. Soc. New York,
65-87(1968).

\bibitem[Rib-Zal]{Rib-Zal} L. Ribes, P. Zalesskii; {\em Prifinite
Groups,} Springer 2000.

\bibitem[Sch]{Sch} M. Schlessinger; {\em Functors of Artin rings,}
Trans. A.M.S.130,208-222(1968).

\bibitem[Tay-Wil]{Tay-Wil} R. Taylor, A. Wiles; {\em
Ring-theoretic properties of certain Hecke algebras,} Ann. of
Math. 142,553-572(1995).

\bibitem[Til]{Til} J. Tilouine; {\em Deformations of Galois representations
and Hecke algebras,} Mehta Research Institute of Allahabad,
Narosa Publishing House 1996.

\bibitem[Tsu]{Tsu} H. Tsungai; {\em On some derivations of Lie algebras
related to Galois representations,} Publ. RIMS. Kyoto Univ.
31,113-134(1995).

\bibitem[Voe]{Voe} V.A. Voevodsky; {\em Galois representations
connected with hyperbolic curves,} Math. USSR Izvestiya
39,1281-1291 (1992).

\bibitem[Wil]{Wil} A. Wiles; {\em Modular elliptic curves and Fermat's last
theorem,} Ann. of Math. 142,443-551(1995).


\end{thebibliography}
\end{document}